\documentclass[11pt,reqno, final]{amsart}

\usepackage{amsmath,amssymb,amsthm,bbm,bm}
\usepackage{color}
\usepackage[colorlinks=true, allcolors=blue]{hyperref}
\usepackage{mathrsfs}
\usepackage{mathtools}

\usepackage[noabbrev,capitalize,nameinlink]{cleveref}
\crefname{equation}{}{}
\usepackage[noadjust]{cite}

\usepackage{fullpage}

\usepackage{graphics}
\usepackage{pifont}
\usepackage[T1]{fontenc}
\usepackage{tikz}
\usetikzlibrary{arrows.meta}
\usepackage{environ}
\usepackage{framed}
\usepackage{url}
\usepackage[linesnumbered,ruled,vlined]{algorithm2e}
\usepackage[noend]{algpseudocode}
\usepackage[labelfont=bf]{caption}
\usepackage{cite}
\usepackage{framed}
\usepackage[framemethod=tikz]{mdframed}
\usepackage{appendix}
\usepackage{graphicx}
\usepackage[textsize=tiny]{todonotes}
\usepackage{tcolorbox}
\usepackage{enumerate}
\allowdisplaybreaks[1]
\usepackage{enumerate}

\crefname{algocf}{Algorithm}{Algorithms}

\crefname{equation}{}{} 
\AtBeginEnvironment{appendices}{\crefalias{section}{appendix}} 

\usepackage[color]{showkeys} 
\colorlet{refkey}{orange!20}
\colorlet{labelkey}{blue!30}

\crefname{algocf}{Algorithm}{Algorithms}

\numberwithin{equation}{section}
\newtheorem{theorem}{Theorem}[section]
\newtheorem{proposition}[theorem]{Proposition}
\newtheorem{lemma}[theorem]{Lemma}

\crefname{claim}{Claim}{Claims}

\newtheorem*{question*}{Question}

\theoremstyle{definition}
\newtheorem{definition}[theorem]{Definition}

\newtheorem*{definition*}{Definition}
\newtheorem{example}[theorem]{Example}

\theoremstyle{remark}
\newtheorem*{remark}{Remark}

\newcommand{\mb}{\mathbb}

\newcommand{\mc}{\mathcal}

\newcommand{\on}{\operatorname}

\renewcommand{\epsilon}{\varepsilon}

\allowdisplaybreaks

\title{Dimension reduction for maximum matchings and the Fastest Mixing Markov Chain}

\author[Jain]{Vishesh Jain}
\address{Department of Statistics, Stanford University}
\email{visheshj@stanford.edu}
\author[Pham]{Huy Tuan Pham}
\address{Department of Mathematics, Stanford University}
\email{huypham@stanford.edu}
\author[Vuong]{Thuy-Duong Vuong}
\address{Department of Computer Science, Stanford University}
\email{tdvuong@stanford.edu}

\begin{document}

\begin{abstract}
Let $G = (V,E)$ be an undirected graph with maximum degree $\Delta$ and vertex conductance $\Psi^*(G)$. We show that there exists a symmetric, stochastic matrix $P$, with off-diagonal entries supported on $E$, whose spectral gap $\gamma^*(P)$ satisfies
\[\Psi^*(G)^{2}/\log\Delta \lesssim \gamma^*(P) \lesssim \Psi^*(G).\]
Our bound is optimal under the Small Set Expansion Hypothesis, and answers a question of Olesker-Taylor and Zanetti, who obtained such a result with $\log\Delta$ replaced by $\log|V|$. 

In order to obtain our result, we show how to embed a negative-type semi-metric $d$ defined on $V$ into a negative-type semi-metric $d'$ supported in $\mb{R}^{O(\log\Delta)}$, such that the (fractional) matching number of the weighted graph $(V,E,d)$ is approximately equal to that of $(V,E,d')$.
\end{abstract}

\maketitle

\section{Introduction}
Let $G = (V,E, w)$ be a simple, undirected graph with non-negative weights $w : E \to \mb{R}_{\geq 0}$. In this note, we consider weights $w$ obtained by restricting $\tilde{w} : V \times V \to \mb{R}_{\geq 0}$, where $\tilde{w}^{1/q}$ embeds isometrically into Euclidean space for some finite $q\geq 1$. In other words, there is a function $f : V \to \mb{R}^{|V|}$ such that for all $u,v \in V$, $\tilde{w}(u,v) = \|f(u) - f(v)\|^{q}$, where $\|\cdot \|$ denotes the usual Euclidean or $\ell_2$ norm. The special case $q=2$ is known as a metric of negative-type.   

Recall that a fractional matching of a weighted graph $G = (V,E,w)$ is a function $h: E \to [0,1]$ such that for all $v \in V$, $\sum_{v\in e}h(e) \leq 1$. A matching corresponds to maps $h: E \to \{0,1\}$. 

Let $M$ (respectively, $M_{\on{frac}}$) denote a maximum weight matching (respectively, fractional matching) of $(G,w)$. Here, the weight of a (fractional) matching $h : E \to [0,1]$ is defined by
\[w(h) := \sum_{e\in E}h(e)w(e).\] The weight of a maximum weight (fractional) matching of $(G,w)$ is called the (fractional) matching number of $(G,w)$. 

It follows from LP duality that the fractional matching number of $(G,w)$ is equal to the value of the following linear program:
\begin{align}
\label{eqn:max-matching-dual}
    \min_{g: V \to \mb{R}_{\geq 0}} &\quad \sum_{u \in V}g(u) \nonumber\\
    \text{s.t. }&\quad g(u) + g(v) \geq w(u,v) \quad \forall \{u,v\} \in E
\end{align}

Motivated by the problem of the Fastest Mixing Markov Chain (see \cref{sec:FMMC}), we consider the following natural dimension reduction question: what is the minimum dimension of a Euclidean space $\mb{R}^{m}$ for which there is a function $f':V \to \mb{R}^{m}$ such that $(G,w') = (V,E,w')$, with $w'(u,v) = \|f'(u) - f'(v)\|^{q}$, has approximately the same (fractional) matching number as $(G,w) = (V,E,w)$? 

A direct application of the classical Johnson-Lindenstrauss lemma shows that $m = O_{q}(\log{n}/\epsilon^2)$ suffices to preserve the (fractional) matching number to within a multiplicative factor of $(1\pm \epsilon)$. Our main result completely removes the dependence on the number of vertices, and shows that $m = O_{q}(\log (\Delta/\epsilon)/\epsilon^2)$ is sufficient, where $\Delta$ denotes the maximum degree of $G = (V,E)$.

\begin{theorem}
\label{thm:matching}
For every $K > 0$, there exists a constant $C_{\ref{thm:FMMC}} = C_{\ref{thm:FMMC}}(K) \geq 1$ for which the following holds. Let $G = (V,E,w)$ with $|V| = n$ and with $w = \tilde{w}|_{E}$, where $\tilde{w}(u,v) = \|f(u) - f(v)\|^{q}$ for some $f : V \to \mb{R}^{n}$ and $q\geq 1$. Let $\Delta$ denote the maximum degree of $G$. Let $M$ (respectively, $M_{\on{frac}}$) denote a maximum weight matching (respectively, fractional matching) of $(G, w)$. 

For $\epsilon \in (0, 1/10)$ and $d \geq C_{\ref{thm:FMMC}}q\log(\Delta/\epsilon q)/\epsilon^{2}$, let $\Pi_{n,d}$ be a $d\times n$ random matrix whose entries are i.i.d.~copies of a centered random variable with variance $1/d$ and sub-Gaussian norm at most $K/\sqrt{d}$. For $u,v \in V$, let $\tilde{w}_{\Pi_{n,d}}(u,v) := \|\Pi_{n,d} f(u) - \Pi_{n,d}f(v)\|^{q}$. Let $M^{\Pi_{n,d}}$ (respectively, $M^{\Pi_{n,d}}_{\on{frac}}$) denote a maximum weight matching (respectively, fractional matching) of $(G, w_{\Pi_{n,d}})$. 

Then, except with probability at most $\exp(- d\epsilon^{2}/\sqrt{C_{\ref{thm:FMMC}}})$ over the realisation $\pi \sim \Pi_{n,d}$, we have
\begin{align*}
    e^{-\epsilon q}w(M) &\leq w_{\pi}(M^{\pi}) \leq e^{\epsilon q}w(M),\\
    e^{-\epsilon q}w(M_{\on{frac}}) &\leq w_{\pi}(M^{\pi}_{\on{frac}}) \leq e^{\epsilon q}w(M_{\on{frac}}).
\end{align*}
\end{theorem}
\begin{remark}
\label{rmk:optimality-matching}
By considering a disjoint union of $n/(\Delta+1)$ stars, each with $(\Delta+1)$ vertices, and weights given by assigning a distinct unit coordinate vector in $\mb{R}^{n}$ to each vertex, we see that $\Omega_{q}(\log(\Delta)/\epsilon^2)$ dimensions are needed for sub-Gaussian dimension reduction schemes of the form considered above. Moreover, if we also demand that $w_{\pi}(E) \geq e^{-\epsilon q}w(E)$ (which holds with high probability in our setting; see \cref{sec:proof-matching}) then the factor of $\log(\Delta)$ in our bound on $d$ cannot be improved (even for possibly non-linear dimension reduction schemes) under the Small Set Expansion Hypothesis (SSE) of Raghavendra and Steurer \cite{raghavendra2010graph}, the reason being that an improved dependence on $\Delta$ in \cref{thm:matching} translates to a corresponding improved dependence on $\Delta$ in \cref{thm:FMMC}, which contradicts SSE (see the remark after \cref{thm:FMMC}).   
\end{remark}

\subsection{Fastest Mixing Markov Chain (FMMC)}
\label{sec:FMMC}

In the Fastest Mixing Markov Chain (FMMC) problem, introduced by Boyd, Diaconis, and Xiao \cite{boyd2004fastest}, we are given a finite, undirected graph $G = (V,E)$ and are asked to design a symmetric $|V|\times |V|$ matrix $P$ such that (i) $P_{ij} \geq 0$ for all $i,j$ and $P_{ij} > 0$ only if $\{i,j\} \in E$ or $i = j$; (ii) $\sum_{j}P_{ij} = 1$ for all $i$; and (iii) denoting the eigenvalues of $P$ in decreasing order by $\lambda_1(P) \geq \lambda_2(P) \geq \dots \geq \lambda_n(P)$, the second largest eigenvalue modulus (SLEM), defined by $$\mu(P) := \max_{i=2,\dots,n}|\lambda_i(P)|,$$ is as small as possible. We denote the set of all symmetric $V \times V$ matrices, satisfying properties (i), (ii), (iii) by $\mc{M}(G)$. 

In words, we are asked to design a discrete-time, time-homogeneous Markov chain with state space $V$, whose transitions are supported on the edges of $E$ and which is reversible with respect to the uniform distribution on $V$, with the smallest possible SLEM among all such Markov chains. Since $\mu(P) \leq 1$, this is equivalent to asking for the largest possible spectral gap, defined by $$\gamma(P) := 1-\mu(P).$$ 

\begin{definition}[Optimal spectral gap]
\label{def:optimal-spectral-gap}
Let $G = (V,E)$ be a finite, undirected graph. With notation as above, we define the optimal spectral gap (among reversible chains supported on $E$) by
\[\gamma^*(G) := \sup\{\gamma(P): P \in \mc{M}(G)\}.\]
\end{definition}

For an extensive discussion of the history of this problem, we refer the reader to \cite[Section~1.4]{olesker2021geometric}, limiting ourselves here to only the most relevant results. All of these results involve the vertex conductance of a graph, whose definition we now recall.

\begin{definition}[Vertex conductance]
\label{def:vertex-conductance}
The vertex conductance $\Psi^*(G)$ of a graph $G = (V,E)$ is defined as
\[\Psi^*(G) := \min_{\emptyset \neq S\subseteq V: |S| \leq |V|/2} \frac{|\{v\notin S : \exists u \in S \text{ s.t. } \{u,v\} \in E\}|}{|S|}\]
\end{definition}

Roch \cite{roch2005bounding} showed that $\gamma^*(G) \lesssim \Psi^*(G)$, thereby identifying the vertex conductance of the host graph as a barrier to designing Markov chains with large spectral gaps. In the other direction, Cheeger's inequality shows that $\gamma^*(G) \gtrsim \Psi^*(G)^{2}/\Delta^{2}$, where $\Delta$ denotes the maximum degree of $G$. Put together, we have
\begin{equation}
    \label{eqn:cheeger}
    \Psi^*(G)^2/\Delta^{2} \lesssim \gamma^*(G) \lesssim \Psi^*(G)
\end{equation}
Simple examples (see \cite{olesker2021geometric}) show that the quadratic dependence on $\Psi^*(G)$ in the lower bound and the linear dependence on $\Psi^*(G)$ in the upper bound cannot be improved in general, leaving open the question of resolving the dependence on $\Delta$ (or removing it altogether). 

Recently, Olesker-Taylor and Zanetti \cite{olesker2021geometric} provided a new lower bound on $\gamma^*(G)$, where the quadratic dependence on $\Delta$ is replaced by logarithmic dependence on $|V|$:

\begin{equation}
\label{eqn:tz}
   \Psi^*(G)^2/\log|V| \lesssim \gamma^*(G) \lesssim \Psi^*(G),
\end{equation}
thereby providing an exponential improvement over \cref{eqn:cheeger} for dense graphs. However, in general, the lower bounds in \cref{eqn:cheeger,eqn:tz} are incomparable, and Olesker-Taylor and Zanetti asked if $\log|V|$ in \cref{eqn:tz} may be replaced by $\log\Delta$, which would provide a common refinement of \cref{eqn:cheeger,eqn:tz} and moreover, would be tight under the Small Set Expansion Hypothesis \cite{raghavendra2010graph} (see the remark after \cref{thm:FMMC}). As an application of \cref{thm:matching}, we answer this question in the affirmative. 

\begin{theorem}
\label{thm:FMMC}
Let $G = (V,E)$ be a graph of maximum degree $\Delta$. Then, $\gamma^*(G)$ satisfies
\[\Psi^*(G)^2/\log\Delta \lesssim \gamma^*(G) \lesssim \Psi^*(G)\]
\end{theorem}
\begin{remark}
\label{rmk:optimality-fmmc}
Work of Louis, Raghavendra, and Vempala \cite{louis2013complexity} shows that under the Small Set Expansion Hypothesis \cite{raghavendra2010graph}, for all sufficiently small $\epsilon > 0$, there is no polynomial-time algorithm to distinguish between $\Psi^*(G) \leq \epsilon$ and $\Psi^*(G) \gtrsim \sqrt{\epsilon \log \Delta}$ for graphs with maximum degree $\Delta$ satisfying $\log\Delta \asymp 1/\epsilon$. As noted in \cite{olesker2021geometric}, a lower bound of the form $\Psi^*(G)^{2}/o(\log(\Delta)) \lesssim \gamma^*(G)$ would allow us to distinguish between the two scenarios (using $\gamma^*(G) \leq \epsilon$ versus $\gamma^*(G) = \omega(\epsilon)$), which would contradict the Small Set Expansion Hypothesis since $\gamma^*(G)$ can be computed in time polynomial in the size of the graph. 
\end{remark}

\begin{example}
A popular use of Markov Chain Monte Carlo methods is to sample from the uniform distribution on an exponentially sized subset $V$ of a product space $\{1,\dots,r\}^{n}$ (where $r \asymp 1$ and $n$ is large) using `local chains'. In our notation, the transitions of the graph are supported on $G = (V,E)$, where $E$ connects vertices within Hamming distance $s$, where $s \asymp 1$. Thus, $\Delta \asymp \log|V| \asymp n$, so that \cref{eqn:cheeger} gives the lower bound $\gamma^*(G) \gtrsim \Psi^*(G)^{2}/n^{2}$, \cref{eqn:tz} gives the lower bound $\gamma^*(G) \gtrsim \Psi^*(G)^{2}/n$, whereas \cref{thm:FMMC} gives the lower bound $\gamma^*(G) \gtrsim \Psi^*(G)^{2}/\log{n}$.   
\end{example}

\subsection{Concurrent work} Shortly after the appearance of our manuscript on the arXiv, we were informed of upcoming work of Kwok, Lau, and Tung \cite{kwok2022cheeger}, which proves \cref{thm:FMMC} using a different proof technique (in particular, \cref{thm:matching} does not appear in \cite{kwok2022cheeger}), as well as a weighted version of \cref{thm:FMMC}, which is not considered in our work.



\section{Proof of \cref{thm:matching}}
\label{sec:proof-matching}
We begin by recalling some preliminary notions related to dimension reduction.

\begin{definition}
\label{defn:good-dimension-reduction}[cf.~\cite[Definition~2.1]{makarychev2019performance}]
Let $\Pi_{n,d}$ denote a distribution on linear maps from $\mb{R}^{n} \to \mb{R}^{d}$. For constants $q\geq 1$ and $c, C > 0$, we say that $\Pi_{n,d}$ is a $(q,c, C)$-good random dimension reduction if for all $1/10 > \epsilon > \frac{C\sqrt{q}}{\sqrt{d}}$, there exist $0 < \delta, \rho < \exp(-c\epsilon^{2}d)$ satisfying the following properties.   
\begin{enumerate}
    \item For all $x,y \in \mb{R}^{n}$, 
    \[\mb{P}_{\pi \sim \Pi_{n,d}}[ \| \pi(x) - \pi(y)\| \notin e^{\pm \epsilon}\|x-y\| ] \leq \delta.\]
    
    \item For all $x,y \in \mb{R}^{n}$, letting $\mc{E}_{x,y}$ be the event that $\|\pi(x) - \pi(y) \| \geq e^{\epsilon}\|x-y\|$, we have
    \[\mb{E}_{\pi \sim \Pi_{n,d}}\left[1_{\mc{E}_{x,y}}\left(\frac{\|\pi(x) - \pi(y)\|^{q}}{\|x-y\|^{q}} - e^{\epsilon q}\right)\right] \leq \rho.\]
\end{enumerate}
\end{definition}

\begin{lemma}[see, e.g., {\cite[Lemma~C.1]{makarychev2019performance}}]
\label{lem:JL}
For every $K_{\ref{lem:JL}} > 0$, there exist constants $c_{\ref{lem:JL}}, C_{\ref{lem:JL}} > 0$ for which the following holds. Let $\Pi_{n,d}$ denote a $d\times n$  random matrix whose entries are i.i.d.~copies of a centered random variable with variance $1/d$ and sub-Gaussian norm at most $K_{\ref{lem:JL}}/\sqrt{d}$. Then, $\Pi_{n,d}$ is a $(q, c_{\ref{lem:JL}}, C_{\ref{lem:JL}})$-good random dimension reduction for all $q\geq 1$.   
\end{lemma}

We are now ready to prove \cref{thm:matching}. Below, we use the notation introduced in the statement of \cref{thm:matching}. For later use, we will prove a slightly stronger statement, which shows that in addition to the conclusion of \cref{thm:matching}, we also have that
\[\sum_{u,v \in V}\tilde{w}_{\pi}(u,v) \geq e^{-\epsilon q}\sum_{u,v \in V}w(u,v).\]

\begin{proof}[Proof of \cref{thm:matching}]
We start by giving a high-level overview of the proof. It is quite immediate to show that the weight of any fixed matching (in particular, a maximum matching) is preserved with high probability under the passage from $w$ to $w_{\pi}$. The challenge, then, is to show that there are no unexpectedly large matchings (in terms of additive error) corresponding to $w_{\pi}$ (note that we cannot use simple union bound due to the extremely large number of matchings). To address this, we consider the (random) set of edges which have unexpectedly high weights under $w_{\pi}$ and show that, with high probability, the total $w_{\pi}$-weight of this random set of edges is at most $\epsilon/\on{poly}(\Delta)$ times the total $w$-weight of all edges and hence, at most $\epsilon/\on{poly}(\Delta)$ times the $w$-matching number of the graph. Thus, while it is possible for the $w_{\pi}$-weight of a matching to be unexpectedly larger (in a multiplicative sense) than its $w$-weight, the magnitude of this increase is much smaller than the $w$-matching number of the graph and so, cannot significantly alter the $w_{\pi}$-matching number. We now proceed to the formal details.     

By \cref{lem:JL}, $\Pi_{n,d}$ is a $(q, c', C')$-good random dimension reduction, where $c', C' > 0$ depend only on $K$. Let $\pi \sim \Pi_{n,d}$, and 
consider the following random subsets of $\binom{V}{2}$, corresponding to heavy edges, light edges, and light pairs. 
\begin{align*}
    \mc{H} &= \{\{u,v\} \in E : w_{\pi}(u,v) \geq e^{\epsilon q}w(u,v)\},\\
    \mc{L}_{1} &= \left\{\{u,v\} \in E : w_{\pi}(u,v) \leq e^{-\epsilon q}w(u,v) \right\},\\
    \mc{L}_{2} &= \left\{\{u,v\} \in \binom{V}{2} : w_{\pi}(u,v) \leq e^{-\epsilon q}w(u,v) \right\}.
\end{align*}
These random subsets naturally give rise to the following quantities corresponding to the excess cost of heavy edges, the total original weight of light edges, and the total original weight of light pairs:
\begin{align*}
\on{Diff}(\mc{H}) := \sum_{\{u,v\} \in \mc{H}}\left(w_{\pi}(u,v) - e^{\epsilon q}w(u,v)\right); \text{  }
\on{Cost}(\mc{L}_1) := \sum_{\{u,v\} \in \mc{L}_1} w(u,v); \text{  } 
\on{Cost}(\mc{L}_2) := \sum_{\{u,v\} \in \mc{L}_2} \tilde{w}(u,v).
\end{align*}
We show that, in expectation, all of these quantities are sufficiently small. First, since any graph of maximum degree $\Delta$ can be written as a disjoint union of at most $(\Delta + 1)$ matchings, it follows that $w(M) \geq \frac{1}{\Delta+1}\sum_{\{u,v\} \in E}w(u,v)$. Then, from \cref{lem:JL,defn:good-dimension-reduction}, it immediately follows that
\begin{align*}
\mb{E}_{\pi}[\on{Diff}(\mc{H})] &\leq \rho \sum_{\{u,v\} \in E}w(u,v) \leq \rho(\Delta+1)w(M),\\
\mb{E}_{\pi}[\on{Cost}(\mc{L}_1)] &\leq \delta \sum_{\{u,v\} \in E}w(u,v) \leq \delta(\Delta+1) w(M), \\
\mb{E}_{\pi}[\on{Cost}(\mc{L}_2)] &\leq \delta \sum_{\{u,v\} \in E}\tilde{w}(u,v) = \delta \tilde{w}\left(\binom{V}{2}\right). 
\end{align*}
Therefore, by Markov's inequality, the union bound, and using $\max{\delta, \rho} \leq \exp(-c'\epsilon^2d)$, we see that except with probability at most $3\exp(-c'\epsilon^2 d/2)$, the following event holds:
\[\mc{G} = \left\{\on{Diff}(\mc{H}) \leq \sqrt{\rho}(\Delta+1)w(M), \on{Cost}(\mc{L}_1) \leq \sqrt{\delta}(\Delta+1)w(M), \on{Cost}(\mc{L}_2) \leq \sqrt{\delta}\tilde{w}\binom{V}{2}\right\}.\]

Let $\pi$ be a realisation of $\Pi_{n,d}$ for which $\mc{G}$ holds. Then, for any (fractional) matching $M'$, we have
\begin{align*}
    w_{\pi}(M') 
    &= \sum_{e}M'(e)w_{\pi}(e) = \sum_{e \notin \mc{H}}M'(e)w_{\pi}(e) + \sum_{e \in \mc{H}}M'(e)e^{\epsilon q} w(e) + \sum_{e \in \mc{H}}M'(e)(w_{\pi}(e) - e^{\epsilon q}w(e))\\
    &\leq e^{\epsilon q}w(M') + \on{Diff}(\mc{H}) \leq e^{\epsilon q}w(M') + \sqrt{\rho}(\Delta + 1)w(M) \\
    &\leq (e^{\epsilon q} + \sqrt{\rho}(\Delta+1))w(M). 
\end{align*}
Moreover, we see that for $M' = M$ or $M' = M_{\on{frac}}$, we have
\begin{align*}
    w_{\pi}(M')
    &= \sum_{e\notin \mc{L}_1}M'(e)w_{\pi}(e) + \sum_{e \in \mc{L}_1}M'(e)w_{e} + \sum_{e\in \mc{L}_1}M'(e)(w_{\pi}(e) - w(e)) \\
    &\geq e^{-\epsilon q}w(M') - \on{Cost}(\mc{L}_1)\\
    &\geq (e^{-\epsilon q} - \sqrt{\delta}(\Delta+1))w(M').
\end{align*}
A similar computation shows that
\begin{align*}
    \tilde{w}_{\pi}\left(\binom{V}{2}\right) \geq (e^{-\epsilon q} - \sqrt{\delta}) \tilde{w}\left(\binom{V}{2}\right).
\end{align*}

Finally, by taking $C$ to be sufficiently large depending on $c', C'$, we can ensure that $$\max\{\sqrt{\delta}, \sqrt{\rho}\}(\Delta+1) \leq \epsilon q/10,$$ so that the desired conclusion follows by rescaling $\epsilon$. \end{proof}

\section{Application to the Fastest Mixing Markov Chain: Proof of \cref{thm:FMMC}}

We will need the following proposition due to Olesker-Taylor and Zanetti. 
\begin{proposition}[\cite{olesker2021geometric}]
\label{prop:fmmc-reduction}
Let $G = (V,E)$ be a simple, undirected graph with $|V| = n$, vertex conductance $\Psi^*(G)$, and optimal spectral gap $\gamma^*(G)$. For $1\leq m \leq n$, let $\lambda_m^*$ denote the optimum value of the following problem:
\begin{align*}
    \min_{f: V \to \mb{R}^{m}, G: V \to \mb{R}_{\geq 0}} &\quad \frac{\sum_{u\in V}g(u)}{\sum_{u \in V}\|f(u)\|^{2}}\\
    \text{s.t. }& \quad \quad \sum_{v\in V}f(v) = 0, \quad \text{and}\\
    & \quad g(u) + g(v) \geq \|f(u) - f(v)\|^{2} \quad \forall \{u,v\} \in E
\end{align*}
Let $m^* = \min\{1 \leq m \leq n : \lambda_n^* \leq \lambda_m^* \leq 2\lambda_n^*\}$. Then, 
\[\Psi^*(G)^2/(\log{m^*})\lesssim \gamma^*(G) \lesssim \Psi^*(G).\]
\end{proposition}

Given the previous proposition, \cref{thm:FMMC} follows from \cref{thm:matching} and LP duality \cref{eqn:max-matching-dual}.

\begin{proof}[Proof of \cref{thm:FMMC}]

By \cref{prop:fmmc-reduction}, and since $\lambda_m^*$ decreases as $m$ increases, it suffices to show that for $d = O(\log\Delta)$, $\lambda_{d}^* \leq 2\lambda^*_n$. To this end, let $F: V \to \mb{R}^{n}$ be such that $\sum_{v\in V}F(v) = 0$, and that the value of the program
\begin{align*}
    \min_{g: V \to \mb{R}_{\geq 0}} & \quad \frac{\sum_{u\in V}g(u)}{\sum_{u \in V}\|F(u)\|^{2}}\\
   \text{s.t.}&\quad  g(u) + g(v) \geq \|F(u) - F(v)\|^{2} \quad \forall \{u,v\} \in E
\end{align*}
is $\lambda_n^*$. 

By LP duality \cref{eqn:max-matching-dual}, $\lambda_n^*$ is the fractional matching number of the graph $G = (V,E,w/S_F)$, where $w(u,v) := \|F(u) - F(v)\|^{2}$ and $S_F := {\sum_{u,v \in V}\|F(u) - F(v)\|^{2}}$. Here, we have used that $\sum_{u,v \in V}\|F(u) - F(v)\|^{2} = \sum_{u \in V}\|F(u)\|^{2}$, which holds since $\sum_{v\in V}F(v) = 0$. 

We now appeal to the slight extension of \cref{thm:matching} (with $K=1$) proved in the previous section. Let $\epsilon = 1/100$, $d = 2C_{\ref{thm:matching}}(1)\log(\Delta/2\epsilon)/\epsilon^{2} = O(\log\Delta)$, and $\pi \in \mc{G}$. Here $\mc{G}$ is the event defined in the previous subsection, where it was also shown to be nonempty. Let $f : V \to \mb{R}^{d}$ be defined by $f(v) = \pi(F(v))$. Then $\sum_{v\in V}f(v) = 0$ and hence, $\lambda_d^*$ is at most 
\begin{align*}
    \min_{g: V \to \mb{R}_{\geq 0}} & \quad \frac{\sum_{u\in V}g(u)}{\sum_{u,v \in V}\|f(u) - f(v)\|^{2}}\\
    \text{s.t.} &\quad g(u) + g(v) \geq \|f(u) - f(v)\|^{2} \quad \forall \{u,v\} \in E.
\end{align*}
Once again, by LP duality \cref{eqn:max-matching-dual}, $\lambda_d^*$ is at most the fractional matching number of the graph $G = (V,E,w'/S'_f)$, where
$w'(u,v) := \|f(u) - f(v)\|^{2}$ and $S'_f := {\sum_{u,v \in V}\|f(u) - f(v)\|^{2}}$.

In the previous section, we showed that for $\pi \in \mc{G}$, $S_f' \geq e^{-2\epsilon}S_F$ and the fractional matching number of $(V,E,w')$ is at most $e^{2\epsilon}$ times the fractional matching number of $(V,E,w)$. By the choice of $\epsilon$, this shows that $\lambda_d^* \leq 2\lambda_n^*$, as desired.

\end{proof}

\bibliographystyle{amsplain0.bst}
\bibliography{main.bib}

\end{document}